\documentclass[12pt]{article}
\usepackage[latin1]{inputenc}
\usepackage{amsmath,amssymb}
\usepackage{latexsym}

\newtheorem{theorem}{Theorem}[section]
\newtheorem{corollary}[theorem]{Corollary}
\newtheorem{lemma}[theorem]{Lemma}
\newtheorem{proposition}[theorem]{Proposition}

\newtheorem{remark}[theorem]{Remark}
\numberwithin{equation}{section}

\def\R{{\mathbb R}}

\def\eps{\varepsilon}

\def\al{\alpha}

\def\phi{\varphi}
\def\qq{\qquad}

\def\proof{{\medskip\noindent {\bf Proof. }}}

\def\qed{{\hfill $\square$ \bigskip}}

\def\sgn{{\mathop {{\rm sgn\, }}}}

\def\square{{\vcenter{\vbox{\hrule height.3pt
        \hbox{\vrule width.3pt height5pt \kern5pt
           \vrule width.3pt}
        \hrule height.3pt}}}}

\def\sA {{\cal A}} \def\sB {{\cal B}} 
  \def\sF {{\cal F}}
  
  \def\sL {{\cal L}}

\def\wh{\widehat}
\def\ol{\overline}
\def\E{{\Bbb E}}
\def\P{{\Bbb P}}
\def\Q{{\Bbb Q}}
\def\norm#1{{\Vert #1 \Vert}}
\def\del{{\partial}}
\def\lam{{\lambda}}
\def\angel#1{{\langle #1 \rangle}}
\def\bea{\begin{align*}}
\def\eea{\end{align*}}
\def\bee{\begin{equation}}
\def\eee{\end{equation}}

\begin{document}

\title{The submartingale problem for a class of degenerate elliptic
operators}

\author{Richard F. Bass\footnote{Research partially supported by NSF grant
DMS-0244737.}  \; and Alexander Lavrentiev}


\maketitle

\begin{abstract}  
\noindent We consider the degenerate elliptic operator
acting on $C^2$ functions on $[0,\infty)^d$: 
\[ \sL f(x)=\sum_{i=1}^d a_i(x) x_i^{\al_i} \frac{\del^2 f}{\del x_i^2}
(x) +\sum_{i=1}^d b_i(x) \frac{\del f}{\del x_i}(x),
\]
where the $a_i$ are continuous functions that are bounded
above and below by  positive constants, the $b_i$ are bounded and
measurable, and the $\al_i\in (0,1)$.
We impose Neumann boundary conditions on the boundary of $[0,\infty)^d$.
There will not be uniqueness for the submartingale problem
corresponding to $\sL$. If  we consider, however,  only those solutions to the
submartingale problem for which the process spends 0 time on the boundary,
 then existence and uniqueness for the submartingale
problem for $\sL$ holds within this class.
Our result is equivalent to establishing weak uniqueness for the system of 
stochastic differential equations
\[
dX_t^i=\sqrt{2a_i(X_t)} (X_t^i)^{\al_i/2}
dW^i_t+b_i(X_t)\, dt +dL_t^{X^i}, \qq X^i_t\geq 0,
\]
where $W_t^i$ are independent Brownian motions and $L^{X_i}_t$ is a local 
time at 0 for $X^i$.

\vskip.2cm
\noindent {\it Keywords:} martingale problem, stochastic differential 
equations, degenerate elliptic operators, speed measure, perturbation,
Bessel process, Littlewood-Paley

\vskip.2cm
\noindent {\it Subject Classification: Primary 60H10; Secondary 60H30}   
\end{abstract}

\section{Introduction}

We consider the degenerate elliptic operator acting on $C^2$ functions on
$[0,\infty)^d$ defined by
\bee\label{Loperator}
\sL f(x)=\sum_{i=1}^d a_i(x) x_i^{\al_i} \frac{\del^2 f}{\del x_i^2}(x)+\sum_{i=1}^d b_i(x)\frac{\del f}{\del x_i}(x).
\eee
We assume here that the $b_i$ are bounded, the $a_i$ are continuous and bounded above
and below by positive constants, and each $\al_i\in (0,1)$. We impose zero Neumann
boundary conditions on $\del(\R^d_+)$, where we write $\R_+=[0,\infty)$.
In this paper we investigate whether there is at most one process corresponding to the operator
$\sL$. 

We formulate this question in terms of a submartingale problem. 
Let $\Omega=C([0,\infty); \R^d_+)$, the continuous functions from $[0,\infty)$ to
$\R^d_+$. Define the canonical process $X$ by $X_t(\omega)=\omega(t)$ and let
$\sF_t$ be the filtration generated by $X$.
Let $x\in \R^d_+$. We say that a probability measure $\P$ on $\Omega$
is a solution to the submartingale problem for $\sL$ started at $x$ if
$\P(X_0=x)=1$ and whenever $f\in C^2(\R_+^d)$ such that for each $i$ we have in
addition $\del f/\del x_i\geq 0$ on 
$\{x=(x_1, \ldots, x_d)\in \R^d_+: x_i=0\}$, 
then
$f(X_t)-f(X_0)-\int_0^t \sL f(X_s) ds$ is a submartingale with respect to $\P$.

If $Y_t$ is a one dimensional process on $[0,\infty)$, by local time at 0 of $Y$
we mean a continuous nondecreasing process $L^Y_t$ such that $L^Y$ increases
only when $Y$ is at 0.
Closely related to the operator $\sL$ is the system of equations
\begin{equation}\label{mainSDE}
dX^i_t= \sqrt{2a_i(X_t)} (X_t^i)^{\al_i/2} dW^i_t + b_i(X_t) dt + dL^{X^i}_t,
\qq i=1, \ldots, d,
\end{equation}
where $X^i_0=x_i$, $X_t^i\geq 0$ for all $t$, $L^{X^i}$ is a local time at 0
for $X^i$, and the $W^i$ are independent one dimensional Brownian motions started
at 0.

We say a weak solution to
(\ref{mainSDE}) exists if there is a probability $\P$ such that (\ref{mainSDE}) holds
and the $W^i$ are independent Brownian motions under $\P$. Weak uniqueness holds
if given any two solutions $(X_j, W_j, \P_j)$, $j=1,2$, the joint law of $(X_1,W_1)$ under $\P_1$
is equal to the joint low of $(X_2, W_2)$ under $\P_2$.

We have assumed that each $\al_i$ is in the interval $(0,1)$, so in fact uniqueness
for the submartingale problem for $\sL$ does {\bf not} hold. This can be seen even
in one dimension: if one looks at the one-dimensional diffusion on natural scale
with speed measure $m(dx)=x^{-\al}$ for $x$ positive and $m$ equal to 0 on $(-\infty,0)$,
one can put an atom of arbitrary finite mass at 0 and obtain different processes. 

If, however, one restricts attention to those solutions to the submartingale
problem $\sL$ for which the process spends zero time at the boundary, then
uniqueness of the submartingale problem does hold. Our main theorem is the
following. Let $\Delta=\del(\R^d_+)$.

\begin{theorem}\label{maintheorem}
Let $x\in \R_+^d$

(a) There exists one and only one solution to the submartingale problem
for $\sL$ started at $x$ that spends zero time in $\Delta$, i.e.,
\[ \int_0^\infty 1_\Delta(X_s) ds=0, \qq \P-\mbox{a.s.}
\]

(b) A weak solution to (\ref{mainSDE}) exists that spends zero time in $\Delta$. Weak uniqueness
holds if we restrict attention to  those weak solutions that spend zero time in $\Delta$.
\end{theorem}

\medskip

Our paper continues the study of degenerate diffusions in the positive orthant begun
in \cite{ABBP} and \cite{BP}. Those papers concerned the operator $\sL$ where all the
$\al_i$ were equal to 1. In some sense, when the $\al_i=1$ is the critical case, in that then the
exact values of the drift coefficients $b_i$ make a large difference to the behavior
of the resulting process. When $\al_i>1$, then either the process never attains the boundary,
or if it starts on the boundary, never leaves, so the problem then becomes a lower dimensional
one. This paper deals with the remaining case.

Although the values of the drift coefficients play less of a role, the results of this paper are
not a subset of those in \cite{ABBP}. In fact, they could not be, because here we need the
additional assumption that the process spends zero time on the boundary in order
to have uniqueness, while no such assumption is needed in \cite{ABBP}. 

If $\al_i<1/4$, one can check that a Girsanov transformation allows one to assume that
the corresponding $b_i$ can be taken to zero. We wanted to allow the full range
of $\al_i$ and drift coefficients, so we did not restrict the values of the $\al_i$
to $(0,1/4)$.
  
As is often the case, uniqueness for a martingale or submartingale problem is often
related to the existence of a solution to a PDE problem. That is also the case here,
but we do not pursue this connection.
Our techniques could also be applied to diffusions on $\R^d$ whose coefficients
decay near the $i^{th}$ axis like $|x_i|^{\al_i}$.
\medskip

Our methods differ substantially from those in \cite{ABBP}. That paper used $L^2$ estimates.
Here we prove an analogue of Krylov's inequality and use Littlewood-Paley theory
to obtain $L^p$ estimates. These are of independent interest, and could be used
to simplify the proof of \cite{ABBP}. In particular, our method
of utilizing  Littlewood-Paley
theory potentially has applications to  many other types of martingale problems.
(The paper \cite{BP} uses $C^\al$ estimates
and is quite different, both in results and in methods.)

The papers \cite{SV2}, \cite{LSz}, and \cite{DI} also consider
diffusions with reflection; the latter two consider pathwise uniqueness. Although some
smoothness of the domain is needed in these papers, the key difference is that 
degeneracies of the type given in (\ref{Loperator}) and (\ref{mainSDE}) are not allowed.

\medskip
After a section on preliminaries, we prove an inequality of Krylov type in Section 3.
Section 4 concerns existence of solutions. Sections 5 and 6 give the required
estimates for the first and second order terms, respectively. The proof of
Theorem \ref{maintheorem} is completed in Section 7. Some of the calculations needed
in Section 5 are deferred to Section 8, which is an appendix.

\section{Preliminaries}

We often write $f_i$ and $f_{ii}$ for the first and second partial derivatives
of $f$ with respect to $x_i$. The Lebsegue measure of a Borel set $B$ will be
written $|B|$. We use $\R_+=[0,\infty)$ and for our state space we will use
$\R_+^d$. We often use $\Delta$ for the boundary of $\R_+^d$. 
More precisely, we set $\Delta_i=\{x=(x_1, \ldots, x_d)\in \R^d_+: x_i=0\}$ and $\Delta=\cup_{i=1}^d \Delta_i$.
We use the letter $c$ with or without
subscripts to denote finite positive constants whose values are unimportant and may change
from place to place.

\medskip

The connection between Theorem \ref{maintheorem} (a) and (b) is the following.

\begin{proposition}\label{weaksol}
There exists a solution to the submartingale problem for $\sL$ started at $x$ if and 
only if there exists a weak solution to (\ref{mainSDE}). The solution to the
submartingale problem will be unique if and only if there is weak uniqueness to
(\ref{mainSDE}). These assertions continue to hold if we restrict attention to
probability measures $\P$ such that $X_t$ spends zero time in $\Delta$, $\P$-a.s.
\end{proposition}

\proof
The proof of this proposition is very similar to the nondegenerate
diffusion case (see \cite{Ba97}, Theorem V.1.1)
and we give only a brief sketch. If $\P$ is a weak solution to (\ref{mainSDE}), then an 
application of Ito's formula shows that $\P$ will be a solution to the submartingale
problem. If $\P$ is a solution to the submartingale problem and 
we take $f(x)=x_i$, then by the definition
of submartingale problem, 
\[
f(X_t)-\int_0^t \sL f(X_s) ds=X_t^i-\int_0^t b_i(X_s) ds
\]
 is a submartingale, and so 
by the Doob-Meyer decomposition can be written as a martingale $M_t^i$  plus an increasing
process $L^i_t$. Similarly to the nondegenerate diffusion case, one can show that $X^i_t-\int_0^t b_i(X_s)\, ds$ is a 
martingale away from the boundary with quadratic variation $\int_0^t a_i(X_s) (X_s^i)^2 ds$.
 This implies that  $L_t^i$ increases only when $X_t^i$ is at 0, and hence is a local time at 0 for $X^i$. 
We then proceed as in the proof of the nondegenerate case.
\qed

Any weak solution to (\ref{mainSDE}) satisfies a tightness estimate. By this we mean
the following.

\begin{proposition}\label{Tight2}
If $M>0, \delta>0, \eta>0$, and $ t_0>0$, there exists $\eps>0$ such that if
$(X,W,\P)$ is any weak solution to (\ref{mainSDE}) and $S_M=\inf\{t : X_t\notin [0,M]^d\}$, then
\[
\P(\sup_{s, t\leq  S_M\land t_0, |t-s|<\eps} |X_t-X_s|>\delta)<\eta.
\]
\end{proposition}

\proof
It suffices to consider each component of $X$ separately. Fix $i$.
If $s<t$ and $X^i_s>\delta/4$, then by
standard estimates 
 we can find $\eps$ such that 
\begin{equation}\label{shorttime}
\P(\sup_{s\leq t\leq s+\eps\leq t_0+\eps} |X^i_t-X^i_s|>\delta/8)<\eta;
\end{equation}
this is because $L^{X^i}$ increases only
when $X^i$ is at 0.
If $X^i_s<\delta/4$, then in order for $X^i$ to be greater than $\delta$
within time $\eps$, there must be times $s'<t'$ with $s<s'<t'<s+\eps$ where $X^i_{s'}=\delta/2$ and
$X^i_{t'}=3\delta/4$. 
But by (\ref{shorttime}),  the probability of this can be made small by taking $\eps$ small.
\qed

When it comes to proving uniqueness of the submartingale problem, it suffices
to consider only solutions defined on the canonical probability space $C([0,\infty); \R_+^d)$.
If $S$ is a stopping time, we let $\Q_S$ be a regular conditional probability $\P(\, \cdot
\circ \theta_S \mid \sF_S)$, where $\theta_S$ is the usual shift operator of Markov
process theory. We denote the corresponding expectation by $\E_{\Q_S}$. Just as
in the nondegenerate case, it is easy to see that $\Q_S$ will be a solution to
the submartingale problem started at $X_S$; cf.\ \cite{Ba97}, Proposition VI.2.1.

\section{Occupation time estimates}

Let 
\bee \label{eai} 
\Delta_i=\{x\in \R^d_+: x_i=0\}, \qquad \Delta=\cup_{i=1}^d \Delta_i.
\eee  
For any process $Z$ and any Borel set $C$ we let
$T_C(Z)=T_C=\inf\{t: Z_t\in C\}$ and $\tau_C(Z)=\tau_C=
\inf\{t: Z_t\notin C\}$. When $C$ is a single point $\{y\}$, we write
instead $T_y$ and $\tau_y$.

We start with an estimate on how long the solution to a one
dimensional SDE can spend near 0.

\begin{theorem}\label{t21} Suppose $x_0\in [0,\infty)$,
$W_t$ is a Brownian motion, $a_t$ and $b_t$ are adapted processes, $c_1^{-1}\leq a_t\leq c_1$ a.s.\ for each $t$,
 and $|b_t|\leq c_1$ a.s.\ for each $t$. Suppose  either (a)
$X$ solves
\bee \label{e21}
 dX_t=a_t\, X_t^{\al/2} dW_t+b_t  \, dt +dL^X_t, \qq X_t\geq 0, \quad X_0=x_0,
\eee
 and $X$ spends zero time at 0

\noindent or (b) for some $\eps>0$
\bee\label{e22}
dX_t=a_t\, (X_t+\eps)^{\al/2} dW_t+b_t \, dt +dL^X_t, \qq X_t\geq 0, \quad X_0=x_0,
\eee
where $L^X_t$ is a continuous nondecreasing process that increases only
when $X_t$ is at 0.
Let $K>0$.
There exists $c_2$ depending only on $K$ and $c_1$
such that
\[
\E \int_0^{T_K(X)} 1_{[0,\eta]}(X_s) ds\leq c_2\eta^{1-\al}.
\]
\end{theorem}

\proof By first performing a nondegenerate time change, we may without
loss of generality suppose that $a_t\equiv 1$.
In case (b), Girsanov's theorem and the fact that
the diffusion coefficient is bounded below away from 0 tells
us that the solution to (\ref{e22}) will
spend zero time at 0. So we can consider both cases at once
if we let $\eps\geq 0$ and specify that $X_t$ spends 0 time
at 0. Let $Y_t=X_t^2$. By Ito's formula,
\bea
dY_t&=2X_t \, dX_t +d\angel{X}_t\\
&=2X_t (X_t+ \eps)^{\al/2} dW_t+ (X_t+\eps)^\al dt
+ 2X_t(X_t+\eps)^{\al/2} b_t\,  dt + 2X_t \, dL^X_t\\
&=2Y_t^{1/2} (Y_t^{1/2}+\eps)^{\al/2} d W_t +
(Y_t^{1/2} +\eps)^\al dt\\
&\qquad  +2 Y_t^{1/2}(Y_t^{1/2}+\eps)^{\al/2}  b_t\, dt,
\end{align*}
where we use the fact that $X_t\,  dL^X_t $ is 0 since $L^X_t$ only increases
when $X_t$ is at 0.

Note $Y_t$ is always nonnegative. By \cite{IW}, Section 6.1,
if $Z_0=Y_0$ and $Z_t$ solves
\[
dZ_t=2Z_t^{1/2} (Z_t^{1/2}+\eps)^{\al/2} d W_t +
(Z_t^{1/2} +\eps)^\al dt -2 c_1 Z_t^{1/2}(Z_t^{1/2}+\eps)^{\al/2} 
dt,
\]
with $Z_t\geq 0$ for all $t$, then $Z_t\leq Y_t$ for all $t$.   

Let $U_t$ be the continuous strong Markov process corresponding to
the operator \[
\tfrac12 (x+\eps)^{\al} f''(x) -c_1  f'(x)
\] 
with reflection at 0;
in the case $\eps=0$ we specify that the process spends 0 time at 0.
Then for some Brownian motion $\wh W$, $U$ will be a weak solution to
\[
dU_t= (U_t+\eps)^{\al/2} d\wh W_t -c_1\sgn(U_t) \, dt+dL^U_t,
\] 
where $L^U_t$ is a continuous nondecreasing process that increases
only when $U_t$ is at 0.
Let $V_t=U_t^2$. 
A calculation using Ito's formula as above shows that $V_t$ solves
the same equation as $Z_t$ except with a different Brownian motion.
Note that the equation defining $Z_t$ locally satisfies the conditions of 
\cite{LeG}, so we have pathwise uniqueness, hence uniqueness in law, and thus
the law of $Z$ and $V$ are the same.  
Therefore
 \bea
\E\int_0^{T_K(X)} 1_{[0,\eta]}(X_s) ds& =\E \int_0^{T_{K^2}(Y)} 1_{[0,\eta^2]}(Y_s) ds\\ &\leq \E\int_0^{T_{K^2}(Z)} 1_{[0, \eta^2]}(Z_s) ds\\ &= \E\int_0^{T_{K^2}(V)} 1_{[0, \eta^2]}(V_s) ds\\
&= \E\int_0^{T_K(U)} 1_{[0,\eta]}(U_s) ds.
\end{align*}

We compute the scale function $s(x)$  for the process
$U_t$ and find that it is determined by
\[
\log s'(x)= \int_0^x \frac{2c_1}{(y+\eps)^\al} dy.
\]
It follows that 
$s(U_t)$ corresponds to the operator 
\[
A_\eps(x) (x+\eps)^\al f''(x)
\]
where $A_\eps(x)$ is a function of $x$ satisfying 
\[
c_3\leq A_\eps(x)\leq c_4, \qquad |x|\leq K,
\]
and $0<c_3<c_4<\infty$ do not depend on $\eps$; furthermore the speed measure
for the process has no atom at 0. Moreover,
we see that $s(\eta)/\eta$ is bounded 
by $c_5$ for $\eta $ small. So 
\[
\E\int_0^{T_K(U)} 1_{[0,\eta]}(U_r) dr
\leq \E\int_0^{T_{s(K)}(s(U))} 1_{[0,s(\eta)]}(s(U_r)) dr.
\]
On the other hand, the term on the right is bounded by
(see \cite{Ba97}, Section IV.3)
\[
c_7\int_{0}^{c_6} 
1_{[0,c_5\eta]}(x) \frac{1}{A_\eps(x) (x+\eps)^\al}
dx
\leq c_8\int_{0}^{c_5\eta} \frac{1}{x^\al} dx,
\]
and this in turn is bounded by
\[ c_9 \eta^{1-\al},\]
independently of $\eps$.
This proves the theorem.
\qed

\begin{corollary}\label{ct21}
Suppose $X_t$ satisfies the hypotheses of Theorem \ref{t21}. Then there exist 
positive $c_1$ and $c_2$
such that for any $\gamma\leq K$, the probability of more than $m$ upcrossings of $[0,\gamma]$
by $X_t$ before time $T_K$ is bounded by 
\[ c_1(1-c_2\gamma)^m. \]
\end{corollary} 

\proof 
For any process $R$ on $\R$, let $T_0^R=0$, $S_i^R=\inf\{t>T_{i-1}^R: R_t=0\}$, $T^R_i
=\{t> S_i^R: R_t=\gamma\}$, $i\geq 1$.
Let $U$ be the process defined in the proof of Theorem \ref{t21}.
Since $s(U)$ is on natural scale,
\bea
\P^\gamma(S_1^U<T_K(U))&=\P^{s(\gamma)}(S_1^{s(U)}<T_{s(K)}(s(U))=1-\frac{s(\gamma)}{s(K)}\\
&\leq 1-c_3\gamma. 
\end{align*}
As in the proof of Theorem \ref{t21}, $X_t$ is stochastically larger than the
process $U_t$. 
Therefore if $X_0=\gamma$,
\[ \P(S^X_1<T_K(X))\leq \P^\gamma(S_1^U<T_K(U))\leq 1-c_3\gamma.
\]
Recall we use $\Q_S$ for a regular conditional probability for $\P(\, \cdot \circ \theta_S\mid \sF_S)$.
Under $\Q_{T_i^X}$ the process $X_t$ satisfies the hypotheses of Theorem \ref{t21}, and so
\bea
\P(T_{i+1}^X\leq T_K(X))&\leq \E[\Q_{T_i^X}(S_1^X<T_K(X)); T_i^X\leq T_K(X)]\\
&\leq (1-c_3\gamma)\P(T_i^X\leq T_K(X)).
\end{align*}
By induction,
\[ \P(T_i^X\leq T_K(X))\leq c_4(1-c_3\gamma)^i, \]
which implies our result.
\qed

\begin{theorem}\label{krylovprop} Let $M>1$ be fixed and $\lam>0$. There exist $p_0$ and $c_1$ depending
on $\al_1, \ldots,\al_d, M, \lam$,  and $d$ such that if $f\in L^{p_0}(\R^{d})$ and $f$ is supported
in $[0,M]^d$, then
\begin{equation}\label{krylov}
\Bigl| \E \int_0^\infty e^{-\lam s} f(X_s)\, ds\Bigr|
\leq  c_1\norm{f}_{p_0}
\end{equation}
for any solution to (\ref{mainSDE}).
\end{theorem}

\proof Let us set $D=[0,2M]^d$ and $E=[0,M]^d$.
Let $A\subset E$ and let $\eps=|A|$. Let $\delta>0$ be a small
positive real to be chosen later and let $F=[\eps^\delta,2M]^d$.

Our first goal is to show that there exists $K_1$ and $\gamma_1$ not depending on $\eps$ such that 
if $A\subset F$, then 
\begin{equation}\label{kk1}
\E \int_0^{\tau_F} 1_A(X_s)\, ds\leq c_2 \eps^{-\delta K_1+\gamma_1}.
\end{equation}
Define 
\[ Y_t^i=[1-\tfrac{\al_i}{2}]^{-1} (X_t^i)^{1-(\al_i/2)}.
\]
We use Ito's formula to see that for $t\leq \tau_F$
\[
dY_t^i= \sqrt 2 a_i(X_t)\, dW_t^i+ \Big[(Y_t^i)^{-\al_i/(2-\al_i)}b_i(X_t) -\frac{\al_i}{4Y_t^i}\Big]\, dt.
\]
Since we are on the set $F$, there is no issue of the degeneracy of $X_t^i$  at 0 causing
problems.
Notice that on $F$ the drift coefficient of $Y^i$ is bounded by $c_{3}\eps^{-K_2\delta}$.
Define the map $\Gamma:F\to \R_+^d$ by
\[ \Gamma(x_1, \ldots, x_d)= ([1-\tfrac{\al_1}{2}]^{-1} (x_1)^{1-(\al_1/2)},
\ldots, [1-\tfrac{\al_d}{2}]^{-1} (x_d)^{1-(\al_d/2)}).
\]
To obtain (\ref{kk1}) it suffices to bound
\begin{equation}\label{kk2}
\E \int_0^{\tau_{\Gamma(F)}} 1_{\Gamma(A)}(Y_s) \,ds. 
\end{equation}
Note
\[ |\Gamma(A)|\leq c_4\prod_{i=1}^d (\eps^\delta)^{-\al_i/2} |A|\leq c_5\eps^{-\delta K_3+1}
\] for some $K_3>0$. 
Let $\ol Y$ be the process whose coefficients agree with those of $Y$ when the process
is in $F$ and satisfy $d\ol Y^i_t=\sqrt 2 |\ol Y_t^i|^{\al_i/2} \, dW^i_t$ 
 outside of $F$. Let $G$ be a ball of
radius $c_6M$ such that contains $F$. 
If we look at the first component of $\ol Y$, 
the time for $|\ol Y^1_t|$ to exceed $2(c_{6}+1)M$ is less
than or equal to  the time $|\ol Y^1_t|$ spends in 
$[0,2(c_{6}+1)M]$ before exceeding $4(c_{6}+1)M$, and this 
latter amount of time has finite expectation
by Theorem \ref{t21}. Therefore
$\E \tau_G(\ol Y)$ is bounded by a constant $c_7$.  We now use Krylov \cite{Krylov}
to obtain the bound
\[ 
\E\int_0^{\tau_G(\ol Y)} 1_{\Gamma(A)}(\ol Y_s)\, ds\leq c_8(1+c_{9}\eps^{-K_2\delta} \E \tau_G(\ol Y))
|\Gamma(A)|^{1/d}.
\]
This inequality follows from a passage to
the limit in equation (4) of \cite{Krylov}. We therefore
have
\begin{align}
\E \int_0^{\tau_F(X)} 1_A(X_s) \, ds
&\leq \E \int_0^{\tau_{\Gamma(F)}(Y)} 1_{\Gamma(A)}(Y_s)\, ds\label{kk3}\\
&\leq \E\int_0^{\tau_{G(\ol Y)}} 1_{\Gamma(A)}(\ol Y_s)\, ds\nonumber\\
&\leq c_{10}(1+\eps^{-K_2\delta}) \eps^{(-\delta K_3+1)/d}. \nonumber
\end{align}
This proves (\ref{kk1}).

Next we will show that if $A\subset E$, then there exists $K_4$, $K_5$, and $\gamma_2$
such that
\begin{equation}\label{kk4}
\E\int_0^{\tau_E} 1_A(X_s)\, ds\leq c_{11} (\eps^{-K_4\delta+\gamma_2}+\eps^{\delta K_5}).
\end{equation}
Write $A=A_1\cup A_2$, where $A_1=A\cap(\cup_{i=1}^d \{0\leq x_i\leq \eps^\delta\})$
and $A_2=A\setminus A_1$. By Theorem \ref{t21},
we know
\begin{equation}\label{kk5}
\E\int_0^{\tau_D} 1_{A_1}(X_s)\, ds\leq \sum_{i=1}^d \E\int_0^{\tau_D} 1_{[0,\eps^\delta]}(X_s^i)\, ds
\leq c_{12} \eps^{\delta K_4}.
\end{equation}
So we need to bound 
\[
\E \int_0^{\tau_D} 1_{A_2}(X_s)\, ds.
\]
Let $T_0=0$, $S_{i}=\inf\{t>T_{i-1}: X_t\in F\}$, and  $T_i=\inf\{t>S_i: X_t\notin [\eps^\delta/2,2M]^d\}$,
$i\geq 1$. Recall that $\Q_{S_i}$ is used for a regular conditional probability. Then
\bea
\E\int_0^{\tau_E} 1_{A_2}(X_s) \, ds
&=\sum_{i=1}^\infty \E\Big[\int_{S_i}^{T_i} 1_{A_2}(X_s); ds; S_i<\tau_E\Big]\\
&=\sum_{i=1}^\infty \E\Big[\E_{\Q_{S_i}} \int_0^{T_1} 1_{A_2}(X_s)\, ds; S_i<\tau_E\Big]\\
&\leq c_{13}\eps^{-K_1\delta +\gamma_1} \sum_{i=1}^\infty \P(S_i<\tau_E).
\end{align*}
Now in order for $S_i$ to be less than $\tau_E$, we must have for some $j\leq d$ at least
$(i-1)/d$ upcrossings of $X^j$ from $\eps^\delta/2$ to $\eps^\delta$ before hitting
the level $M$. We know by Corollary \ref{ct21} that the probability of this happening is
less than $c_{14} (1-c_{15}\eps^\delta)^{(i-1)/d}$. 
Therefore 
\[
\sum_{i=1}^\infty \P(S_i<\tau_E) \leq c_{16} \eps^{-K_5\delta}.
\]
Combining with (\ref{kk5}), we
have
\begin{equation}\label{kk6}
\E\int_0^{\tau_E} 1_A(X_s)\, ds\leq c_{17}(\eps^{\delta K_4}+\eps^{-K_1\delta+\gamma_1-K_5\delta}),
\end{equation}
which yields (\ref{kk4}).

The next step is to show that there exists $K_6$ and $\gamma_3$ such that
if $A\subset E$, then
\begin{equation}\label{kk7}
\E\int_0^\infty e^{-\lam s} 1_A(X_s)\, ds \leq c_{18} (\eps^{-K_6\delta+\gamma_3}
+\eps^{\delta K_4}).
\end{equation}
Let $V_0=0$, $U_i=\inf\{t>V_{i-1}: X_t\in E\}$, $V_i=\inf\{t>U_i: X_t\notin D\}$, $i\geq 1$. 
We have for $x\notin [0,2M)^d$ that there exists $c_{19}$ and $c_{20}$ such that
$\P^x(U_1>c_{19})>c_{20}$. This is because if $x\notin [0,2M)^d$, then at least one coordinate
of $X$ is greater than or equal to $2M$, and in this range this coordinate is a diffusion 
whose diffusion coefficients are bounded above and below and whose drift coefficient
is bounded above; therefore it cannot move a distance $M$ too quickly. We conclude
\[
\E e^{-\lam U_1}\leq (1-c_{20})+c_{20}e^{-\lam c_{19}}:=\rho.
\]
Note that $\rho<1$. We then have 
\bea
\E e^{-\lam U_i}&=\E\Big[e^{-\lam V_{i-1}}\E\Big[e^{-\lam  U_i}\mid \sF_{V_{i-1}}\Big]\, \Big]\\
&\leq \E\Big[ e^{-\lam U_{i-1}}\E_{\Q_{V_{i-1}}} e^{-\lam U_1}\Big]
\leq \rho \E e^{-\lam U_{i-1}}.
\end{align*}
So by induction
\[
\E e^{-\lam U_i}\leq \rho^{i-1}.
\]
Then
\bea
\E\int_0^\infty e^{-\lam s} 1_A(X_s)\, ds
&=\sum_{i=1}^\infty \E \int_{U_i}^{V_i} e^{-\lam s} 1_A(X_s)\, ds\\
&=\sum_{i=1}^\infty e^{-\lam U_i} \E_{\Q_{U_i}}\int_0^{V_1} e^{-\lam s} 1_A(X_s)\, ds\\
&\leq c_{21}( \eps^{-K_4\delta+\gamma_2} +\eps^{K_4\delta})\sum_{i=1}^\infty  \E e^{-\lam U_i}\\
&\leq c_{22}(\eps^{-K_4\delta+\gamma_2}+\eps^{K_4\delta}),
\end{align*}
which proves (\ref{kk7}).

We now choose $\delta=\gamma_3/(2K_6)$ and we obtain
\begin{equation}\label{kk8}
\E\int_0^\infty e^{-\lam s} 1_A(X_s) \, ds \leq c_{23} \eps^{\gamma_4}.
\end{equation}

Let $p_0=2/\gamma_4$.
To complete the proof, suppose $f\in L^{p_0}(\R_+^d)$ with support in $[0,M]^d$.
By multiplying by a constant, it suffices to consider the
case where  $\norm{f}_{p_0}=1$.
Without loss of generality we may also suppose $f\geq 0$.
Let $A_n=\{x: f(x)\geq 2^n\}$. Then 
\[ 
|A_n|\leq \frac{\norm{f}_{p_0}^{p_0}}{(2^n)^{p_0}}=2^{-n{p_0}}.
\]
We then have
\[
\E\int_0^\infty 1_{A_n}(X_s)\, ds\leq c_{24} 2^{-n{p_0}\gamma_4}.
\]
Thus 
\bea
\E\int_0^\infty e^{-\lam s} f(X_s)\, ds&\leq 1+\sum_{n=0}^\infty 2^{n+1}\E\int_0^\infty
e^{-\lam s} 1_{A_n}(X_s) \,ds\\
&\leq 1+\sum_{n=0}^\infty 2^{n+1}c_{24} 2^{-n{p_0}\gamma_4}\leq c_{25}<\infty.
\end{align*}
Since $\norm{f}_{p_0}=1$, the proof is complete.
\qed

\section{Existence}

In this section we prove existence of a solution to the submartingale problem.
 for the operator $\sL$ defined in (\ref{Loperator}) 
with Neumann boundary conditions on $\Delta$.

There are two complications that are not present in the usual
case: we need to show that our solution spends zero time on
the set $\Delta$ defined in (\ref{eai}); and unless the $\al_i$ are small,
 we cannot use the Girsanov transformation to 
reduce to  the case of zero drift.

Let $\sL^\eps$ be the operator defined by
\[
\sL^\eps f(x)=\sum_{i=1}^d a_i(x) (x_i+\eps)^{\al_i} f_{ii}(x)
+\sum_{i=1}^d b_i(x) f_i(x),
\]
again with reflecting boundary conditions on $\Delta=\del(\R^d_+)$.
The diffusion coefficients are uniformly positive definite and continuous
and are of at most linear growth,
so there exists a unique solution to the submartingale problem for
$\sL^\eps$ started from $x_0$ for every $x_0$; let us denote it
$\P_\eps$. (We reflect the coefficients over the coordinate axes, construct
the solution to the corresponding martingale problem on $\R^d$, and then
look at the law of $(|X^1|, \ldots, |X^d|)$.)

Using Proposition \ref{Tight2} it is standard  that the $\P_\eps$ are a tight sequence of 
probability measures on $C([0,\infty); \R^d_+)$ and there must
exist a subsequence $\eps_j$ such that $\P_{\eps_j}$ converges
weakly. Denote the limit measure by $\P$ and the corresponding 
expectation by $\E$. It is obvious
that $\P(X_0=x_0)=1$.

\begin{proposition}\label{zerotime}
Under $\P$ the process spends zero time in $\Delta$, i.e.,
\[
 \int_0^\infty 1_\Delta(X_s) ds=0, \qquad \P-\mbox{a.s.}
\]
\end{proposition}

\proof
Under $\P_\eps$, the $i$th component
of $X_t$ will satisfy an SDE of the form
\[
dX_t^i= \sqrt {2 a_i(X_t) } (X_t^i+\eps)^{\al_i/2}  dW^i_t
+ b_i(X_t )\, dt +dL^{X^i}_t,
\]
where $L^{X^i}$ is a local time at 0.
Applying Theorem \ref{t21} the amount of time
$X_t^i$ spends in $[0,\eta]$ before exceeding $K$ 
 under $\P_\eps$ is bounded by $c_1 \eta^{1-\al_i}$, where $c_1$
may depend on $K$ but not $\eps$. Taking a limit
\[
\E_\P \int_0^{\tau_K(X^i)}1_{[0,\eta]} (X_s^i) ds
\leq c_1 \eta^{1-\al_i}.
\]
Letting $\eta\to 0$, we have
\[
\E_\P\int_0^{\tau_K(X^i)} 1_{\Delta_i}(X_s) \, ds=0.
\]
Since $K$ is arbitrary, and using this argument for each $i=1, \ldots, d$,
we obtain
\bee\label{e31}
\E_\P \int_0^\infty 1_\Delta(X_s) ds=0,
\eee
and hence the amount of time spent in $\Delta$ is 0 almost surely.
\qed

\begin{proposition}\label{solution}
$\P$ is a solution to the submartingale problem for $\sL$ started
at $x_0$.
\end{proposition}

\proof
To prove that $\P$ is a solution to the submartingale problem
for $\sL$ started at $x_0$, we need to prove that
$M^f_t=f(X_t)-f(X_0)-\int_0^t \sL f(X_s) ds$ is a submartingale
for every $f\in C^2$ such that $f_i\geq 0$ on $\Delta_i$. To do that, it suffices to show
$M^f_t$ is a submartingale for every such $f\in C^2$ that in addition has compact support.
Take such an $f$. It will suffice to show
\bee \label{e32}
\E[ M_t^f Y]\ge \E[M_s^f Y]
\eee
whenever $s<t$,
$n>0$, $r_1\leq r_2\leq \cdots \leq r_n\leq s$, and
$Y=\prod_{i=1}^n g_i(X_{r_i})$ where $g_1, \ldots, g_n$ are bounded
continuous functions with compact support. We know
\bea
\E_\eps&\Big[ \Bigl\{f(X_t)-f(X_0)-\int_0^t \sL^\eps f(X_u) du\Bigr\} Y\Big]\\
&\qquad \geq \E_\eps\Big[ \Bigl\{f(X_s)-f(X_0)-\int_0^s \sL^\eps f(X_u) du\Bigr\} Y\Big]
\end{align*}
since $f(X_t)-f(X_0)-\int_0^t \sL^\eps f(X_u) \, du$ is a submartingale under
$\P_\eps$.
Since $f(X_s)Y$, $f(X_t)Y$, and $f(X_0)Y$ are each continuous functionals of the path,
then the expectations under $\P_{\eps_j}$ converge to the expectations under $\P$, and 
we will be done if we prove
\bee\label{e33}
\E_{\eps_j}\Big[Y\int_0^t \sL^{\eps_j} f(X_u) du\Big]\to
\E\Big[Y\int_0^t \sL f(X_u) du \Big]
\eee
and also with $t$ replaced by $s$. 
We will do only the $t$ case; the $s$ case is almost identical.

Since $f$ is $C^2$ with compact support, then $\sL^\eps f(x)\to \sL f(x)$ uniformly.
So
it suffices to prove
\bee\label{e34}
\E_{\eps_j}\Big[Y\int_0^t \sL f(X_u) du\Big]\to
\E\Big[Y\int_0^t \sL f(X_u) du \Big].
\eee
Let $F(x)=\sum_{i=1}^d b_i(x) f_i(x).$
Since the $a_i$ are continuous, 
then the $\P^{\eps_j}$ expectations of
\[
Y\int_0^t \sum_{i} a_i(X_u) du
\] converge to the expectation under $\P$, and it suffices to prove
\bee\label{e35}
\E_{\eps_j}\Big[Y\int_0^t F(X_u) du\Big]\to
\E\Big[Y\int_0^t F(X_u) du \Big].
\eee

Let $\delta>0$. Pick $K$ large so that 
\[ 
\P_\eps(T_K<t)<\delta.
\]
By tightness, this can be done uniformly in $\eps$ and the same
inequality holds with $\P$ in place of $\P_\eps$. 
So since $\delta$ is arbitrary, it suffices to show
\bee\label{e36}
\E_{\eps_j}\Big[Y\int_0^t  G(X_u) du\Big]\to
\E\Big[Y\int_0^t G(X_u) du \Big],
\eee
where $G(x)=F(x)$ for $|x|\leq K$ and 0 otherwise.
Let $\gamma>0$.
We can find a continuous bounded function $H$ that is 0 outside of
$B(0,K)$, that is equal to
$G$ except on a set of Lebesgue measure less than $\gamma$, and
where $\norm{H}_\infty=\norm{G}_\infty$.
Since $H$ is continuous,
\bee \label{e37}
\E_{\eps_j}\Big[Y \int_0^t H(X_u) du\Big] \to \E\Big[Y\int_0^t H(X_u) du\Big].
\eee
But by Theorem \ref{krylovprop} we have
\[
\E_{\eps} \int_0^{T_K} |H(X_u)-G(X_u)| du
\leq c_2(K,\gamma),
\]
uniformly in $\eps$, 
which can be made less than $\delta$ if we take $\gamma$ small
enough.
Passing to the limit along $\eps_j$,
we have the same result when $\E_\eps$ is replaced by $\E$.
This, (\ref{e37}), and the facts that $\P_\eps(T_K<t)<\delta$ and 
$\P(T_K<t)<\delta $ suffice to establish (\ref{e36}).
\qed

\section{First order estimates}

We first consider the continuous strong Markov process
$Z_t$ on $[0,\infty)$ associated with the operator
\[
\sA_Z f(x)=x^\al f''(x).
\]
Here $\al\in (0,1)$ and we impose reflecting boundary conditions at 0.
More precisely, we  have a process on natural scale whose speed measure
has no atom at 0 and does not charge $(-\infty,0]$. 

Let 
\[
b=1-\frac{\al}{2},
\]
and note that $b\in (\frac12,1)$.
If we set 
\[Y_t=\frac{1}{b\sqrt 2 }Z_t^b,
\]
 a straightforward calculation shows that
$Y$ is a continuous strong Markov process on $[0,\infty)$ associated to the operator
\[
\sA_Y f(x)= \tfrac12 f''(x) +\frac{b-1}{2b x} f'(x)
\]
with reflection at 0, i.e., a Bessel process of order $\delta=\frac{b-1}{b}+1$
with reflection at 0. By \cite{RY} the transition densities of $Y$ (with respect
to Lebesgue measure) are given by
\[
p_Y(t,x,y)=  \Big(\frac{y}{x}\Big)^\nu
\frac{y}{t}e^{-(x^2+y^2)/2t} I_\nu(xy/t),
\]
where
$\nu=\frac{\delta}{2}-1=-\frac{1}{2b}$ and $I_\nu$ is the standard modified
Bessel function.

A change of variables then gives
\bee\label{ptdef}
p_Z(t,x,y)=\frac{c_1}{t}y^{2b-\frac32} e^{-c_2y^{2b}/2t} x^{\frac12} e^{-c_2x^{2b}/2t} I_\nu(c_2x^by^b/t)
\eee
and we have the scaling relationship
\bee \label{ptscale}
p_Z(t,x,y)=t^{-1/2b}p_Z(1,xt^{-1/2b},yt^{-1/2b}).
\eee

We will need the following lemma, the proof of which is given in the appendix.

\begin{lemma} \label{l41} There exists a constant $c_1$ such that
\begin{align}
\sup_x \int_0^\infty \Bigl|\frac{\del }{\del x}p_Z(t,x,y)\Bigr| dy&\leq c_1t^{-\frac{1}{2b}}; \label{ap23}\\
\sup_y y^\al \int_0^\infty \Bigl|\frac{\del}{\del x} p_Z(t,x,y) \Bigr| x^{-\al} dx&
\leq c_1t^{-\frac{1}{2b}}. \label{ap25}
\end{align}
\end{lemma}

Let $P_t$ be the semigroup for $Z_t$, i.e.,  $P_tf(x)=\E^x f(Z_t)$. Let
$\mu(dx)=x^{-\al} dx$ and we consider the space 
$L^2(\R_+, \mu)$.  We use  the above estimates to prove

\begin{proposition}\label{p44}
Suppose $p\in (1,\infty)$. There exists a constant $c_1$ depending only on $p$  such that
\[
\norm{(P_t f)'}_p\leq c_1 t^{-1/2b} \norm{f}_p, \qquad f\in L^2(\R_+,\mu).
\]
\end{proposition}

\proof Fix $t>0$ and write $K(x,y)$ for $|\frac{\del}{\del x} p(t,x,y)|$. 
Let $q$ be the conjugate exponent to $p$. Then by
Lemma \ref{l41}  we have
\bea
\norm{(P_tf)'}_p^p&\leq \int_{\R_+} \Big[ \int_{\R_+} K(x,y) |f(y)|\, dy\Big]^p x^{-\al} dx\\
&= \int_{\R_+} \Big[ \int_{\R_+} K(x,y)^{1/q} K(x,y)^{1/p} |f(y)|\, dy\Big]^p x^{-\al} dx\\
&\leq \int_{\R_+} \Big[\int_{\R_+} K(x,y) \, dy\Big]^{p/q} \Big[\int_{\R_+} K(x,y) |f(y)|^p \, dy\Big] x^{-\al}dx\\
&\leq c_2t^{-(1/2b)(p/q)} \int_{\R_+}\int_{\R_+} K(x,y) x^{-\al}\, dx\,|f(y)|^p\, dy\\ 
&\leq c_3t^{-(1/2b)((p/q)+1)} \int_{\R_+} |f(y)|^p y^{-\al}\, dy.
\end{align*}
Now take $p$-th roots of both sides.
\qed

\bigskip
Now let us turn to the $d$-dimensional case. We suppose $Z_t^i$
is the process on $\R$ corresponding 
to the operator 
\[
\sA_i f(x)=x^{\al_i} f''(x),
\]
with the speed measure having no atom at 0 and not charging 
$(-\infty,0)$ and $\al_i\in (0,1)$. Set
$b_i=1-\frac{\alpha_i}{2}$. We let $p_i(t,x,y)$
denote the transition densities of $Z_t^i$, $P^i_t$ the corresponding
semigroups, let $Z_t=(Z_t^1, \ldots, Z_t^d)$, and
let $p(t,x,y)=\prod_{i=1}^d p_i(t,x_i,y_i)$ be the transition densities for
$Z$ when $x=(x_1, \ldots, x_d)$, $y=(y_1, \ldots, y_d)$. Let $P_t$ now denote
the semigroup for $Z$. Let
$\mu_i$ be the measure on $\R_+$ whose Radon-Nikodym derivative with respect to
$d$-dimensional Lebesgue measure is $x_i^{-\al_i}$ and let $\mu$ be the measure
on $\R_+^d$ given by
\begin{equation} \label{defmu}
\mu(dx)=\prod_{i=1}^d \mu_i(dx_i).         
\end{equation}

We have the analogue of Proposition \ref{p44}.

\begin{proposition}\label{p45}
There exists a constant $c_1$ such that
for each $i$
\[
\Bigl\Vert{\frac{\del(P_t f)}{\del x_i}}\Bigr\Vert_p\leq c_1 t^{-1/2b_i} \norm{f}_p, \qquad f\in L^p(\R_+^d,\mu).
\]
\end{proposition}

\proof We will prove this in the case $i=1$, the case for other $i$'s being exactly
similar. Let $\ol x=(x_2, \ldots, x_d)$. Let $f\in L^2(\R_+^d, \mu)$ and set
\[
F(x_1;\ol x)=
\int\cdots\int \prod_{j=2}^d p_j(t,x_j,y_j) f(y_1,y_2, \ldots y_d)\, dy_2\, \cdots\, dy_d.
\]
Then
\[
\Bigl|\frac{\del}{\del x_1} P_tf(x)\Bigr|=\Bigl|\frac{\del}{\del x_1} P^1_t F(x_1; \ol x)\Bigr|, 
\] and so by Proposition \ref{p44} we have
\[
\int \Bigl|\frac{\del P_tf(x)}{\del x_1}\Bigr|^p\mu_1(dx_1)\leq
c_2t^{-p/2b_1} \int |F(x_1; \ol x)|^p \mu_1(dx_1).
\]
If we integrate both sides with respect to $\ol \mu(dx_2\cdots dx_d)
=\prod_{j=2}^d \mu_j(dx_j)$, we will
have our result provided we show
\bee\label{e44}
\int |F(x_1; \ol x)|^p \mu(dx)\leq \int |f|^p \mu(dx).
\eee

To prove (\ref{e44}) let $\ol P_t$ be the semigroup corresponding to $(Z_t^2, \ldots, Z_t^d)$.
It is easy to check that $\sum_{j=2}^d \sA_j$ is self-adjoint with respect to
the measure $\ol \mu$. Therefore, using Jensen's inequality,
\[
\norm{\ol P_t g}_{L^p(\ol \mu)}\leq \norm{g}_{L^p(\ol \mu)}, \qquad g:\R_+^{d-1}\to \R,
\]
or
\bee \label{e45}
\int |\ol P_t g(x)|^p \ol \mu(dx)\leq \int |g(x)|^p \ol \mu(dx).
\eee
We hold $x_1$ fixed and apply this to $g(\ol x; x_1)=f(x_1, \ldots, x_d)$.
Note $\ol P_t g(\cdot; x_1)=F(x_1; \ol x)$.
So applying (\ref{e45}) to this $g$, we have
\[
\int |F(x_1; \ol x)|^p \ol \mu(dx)\leq \int |f(x_1, \ldots, x_d)|^p \ol \mu(dx).
\]
(\ref{e44}) follows by integrating both sides of this equation
 with respect to $\mu_1(dx_1)$.
\qed

Our main result of this section is the following.
Let 
\[
R_\lam f=\int_0^\infty e^{-\lam t} P_t f\, dt.
\]

\begin{theorem}\label{t46}
There exists $c_1$ such that
for each $i$
\[
\norm{\del (R_\lam f)/\del x_i}_p\leq c_1\lam^{\frac{1}{2b_i}-1}\norm{f}_p.
\]
\end{theorem}

\proof Since $-1/(2b_i)>-1$, the result follows from Proposition \ref{p45},
dominated convergence, and
Minkowski's inequality for integrals:
\begin{align*}
\Bigl\| \frac{\del}{\del x_i} \int_0^\infty e^{-\lam t} P_t f\, dt
\Bigr\|_p&=\Bigl\| \int_0^\infty  e^{-\lam t} \frac{\del P_t f}{\del x_i} dt
\Bigr\|_p \\
&\leq \int_0^\infty e^{-\lam t} c_2 t^{-1/(2b_i)} dt \norm{f}_p.
\end{align*}
\qed

\begin{remark}\label{rem5}               
{\rm Only very minor changes are needed to get the same conclusion
as in Theorem \ref{t46} if we instead set $R_\lam$ to be the resolvent
for the operator $\sum_{i=1}^d a_i \sA_i$, where the $a_i$ are strictly
positive finite constants.}
\end{remark}

\section{Second order estimates}

Let 
\[ 
\sA_i f(x)=|x_i|^{\al_i} f_{ii}(x), \qq x\in \R,
\]
 and let $P_t^i$ be the semigroup corresponding
to the process $Y^i_t$ associated with  $\sA_i$ that spends zero time at 0.
We let $P_t$ be the semigroup corresponding to the process $Y_t=(Y^1_t, \ldots, Y^d_t)$,
where the $Y^i_t$ are independent. The independence implies that if $f(x)=\prod_{i=1}^d f^{(i)}(x_i)$
and $x=(x_1, \ldots, x_d)$,
then $P_t f(x)=\prod_{i=1}^d P_t^i f^{(i)}(x_i)$.

We let $U_t$ be the Poisson semigroup defined in terms of $P_t$:
\[ 
U_t =\int_0^\infty   \Big(\frac{t}{2\sqrt  \pi} e^{-t^2/4s} s^{-3/2}\Big)\, P_s\, ds;
\]
see  \cite{Meyer1}, p.~127. 

The semigroup $P_t^i$ is self-adjoint on $(\R, \mu_i)$, where
$\mu_i(dx)=|x_i|^{-\al_i} \, dx_i$. We let $\mu(dx)=\prod_{i=1}^d \mu_i(dx_i)$ be the
product  measure
on $\R^d$. 

Define
\[ R_0f(x)=\int_0^\infty P_t f(x)\, dt.
\]

\begin{lemma}\label{lpident}
For $f, h$ be bounded on $\R^d$  with compact support 
we have the identity
\begin{equation}\label{lpidenteq}
\int (\sA_i R_0 f(x)) h(x)\, \mu(dx)=
\int \int_0^\infty y(\sA_i U_{y/2}f(x))(U_{y/2}h(x))\, dy\, \mu(dx).
\end{equation}
\end{lemma}

\proof Using the spectral theorem, there exists  (see \cite{Meyer2}) a spectral representation 
\[ P_t^i = \int_0^\infty e^{-\lam_i t} dE^i_{\lam_i}, \qq i=1, \ldots, d. \]
Write $s(\lam)=\sum_{i=1}^d \lam_i$ if $\lam=(\lam_1, \ldots, \lam_d)$. If $f(x)=\prod_{i=1}^d f^{(i)}(x_i)$,
then
\[ P_t f=\int_0^\infty \cdots \int_0^\infty e^{-ts(\lam)} dE^1_{\lam_1}(f^{(1)}) \cdots
dE^d_{\lam_d}(f^{(d)}), \]
and so
\[ R_0f=\int_0^\infty\cdots \int_0^\infty \frac{1}{s(\lam)}dE^1_{\lam_1}(f^{(1)}) \cdots
dE^d_{\lam_d}(f^{(d)}). \]
We have (\cite{Meyer1}, p.~127)
\[ U_t f= \int_0^\infty\cdots \int_0^\infty e^{-t\sqrt{s(\lam)}}dE^1_{\lam_1}(f^{(1)}) \cdots
dE^d_{\lam_d}(f^{(d)}). \]
Note also 
\[ \sA_i =\int_0^\infty \lam_i dE^i_{\lam_i}. \]
Therefore, if $h(x)=\prod_{i=1}^d h^{(i)}(x_i)$, the left hand side of
(\ref{lpidenteq}) is
\begin{equation}\label{lp1}
\int_0^\infty\cdots \int_0^\infty \frac{\lam_i}{s(\lam)} d(E^1_{\lam_1} (f^{(1)}),E^1_{\lam_1}(h^{(1)}))
\cdots d(E^d_{\lam_d} (f^{(d)}),E^d_{\lam_d}(h^{(d)})).
\end{equation}
We use here $(\, \cdot, \cdot\,)$ for the inner product in $L^2(\mu)$.

Similarly, the right hand side of (\ref{lpidenteq}) is
\begin{align*}
\int_0^\infty\,\Big[ \int_0^\infty\cdots &\int_0^\infty y\, \lam_i e^{-(y/2)\sqrt{s(\lam)}}
e^{-(y/2)\sqrt{s(\lam)}} d(E^1_{\lam_1} (f^{(1)}),E^1_{\lam_1}(h^{(1)}))\\
&\cdots d(E^d_{\lam_d} (f^{(d)}),E^d_{\lam_d}(h^{(d)}))\Big]\, dy\\
&=\int_0^\infty\,\Big[ \int_0^\infty\cdots \int_0^\infty y\, \lam_i e^{-y\sqrt{s(\lam)}} 
d(E^1_{\lam_1} (f^{(1)}),E^1_{\lam_1}(h^{(1)}))\\
&\qq \cdots d(E^d_{\lam_d} (f^{(d)}),E^d_{\lam_d}(h^{(d)}))\Big]\, dy.
\end{align*} 
Since 
\[ \int_0^\infty ye^{-Ky}\, dy=\frac{1}{K^2}, \]
this is equal to
(\ref{lp1}). Linear combinations of functions of the form $\prod_{i=1}^d f^{(i)}(x_i)$
 are dense
in $L^2(\mu)$, and an approximation argument completes the proof.
\qed

Define for $(x,y)\in \R^d\times [0,\infty)$
\begin{equation}\label{lpgfunc}
G(f)(x)=\Bigg(\int_0^\infty y\Big(\Bigl|\frac{\del U_yf}{\del y}(x,y)\Bigr|^2
+\sum_{i=1}^d |x_i|^{\al_i}\Bigl|\frac{\del U_yf}{\del x_i}(x,y)\Bigr|^2\Big)\, dy\Bigg)^{1/2}.
\end{equation}
The main result of \cite{Meyer3} implies
that if $1<p<\infty$, then there exists $c_p$ such that
\begin{equation}\label{lpineq}
\norm{G(f)}_p\leq c_p \norm{f}_p,
\end{equation}
where the norm is the $L^p(\mu)$ norm.

The main theorem of this section is the following.

\begin{theorem}\label{lptheo}
Let $1<p<\infty$. There exists a constant $c_1$ depending only on $p$ such that
\[ \norm{\sA_i R_0 f}_p\leq c_1 \norm{f}_p.
\]
\end{theorem}

\proof Let $f$ and $h$ be bounded with compact support. Using Lemma \ref{lpident},
integration by parts, Cauchy-Schwarz, 
and
H\"older's inequality, we have
\begin{align*}
\int(\sA_i &R_0 f(x)) h(x)\, \mu(dx)\\
& =\int \int_0^\infty y(\sA_i U_{y/2}f(x))(U_{y/2}h(x))\, dy\, \mu(dx)\\
&=\int\int_0^\infty y |x_i|^{\al_i} \frac{\del U_{y/2}f}{\del x_i}(x)
\frac{\del U_{y/2}h}{\del x_i}(x)\, dy\, \mu(dx)\\
&\leq \int \Big(\int_0^\infty y |x_i|^{\al_i}\Bigl|\frac{\del U_{y/2}f}{\del x_i}(x)\Bigr|^2\, dy\Big)^{1/2}\\
&\qq \times \Big(\int_0^\infty y |x_i|^{\al_i}\Bigl|\frac{\del U_{y/2}h}{\del x_i}(x)\Bigr|^2\, dy\Big)^{1/2}\, \mu(dx)\\
&\leq \norm{G(f)}_p \norm{G(h)}_q,
\end{align*}
where $q$ is the conjugate exponent to $p$.
By (\ref{lpineq}) this in turn is bounded by
\[ c_2\norm{f}_p\norm{h}_q. \]
If we now take the supremum over all such $h$ for which $\norm{h}_q\leq 1$, we obtain
\[ \norm{\sA_i R_0f}_p\leq c_3\norm{f}_p \]
for $f$ bounded with compact support. An approximation argument allows us to extend
this inequality to all $f\in L^p$.
\qed

\begin{corollary}\label{lpcor}
Let $\lam>0$.
Let $1<p<\infty$. There exists a constant $c_1$ depending only on $p$ such that
\[ \norm{\sA_i R_\lam f}_p\leq c_1 \norm{f}_p.
\]
\end{corollary}

\proof 
If $f\in L^p$, then $f-\lam R_\lam f$ is also in $L^p$ with $\norm{f-\lam R_\lam f}_p\leq 2\norm{f}_p$.
Our result now follows from Theorem \ref{lptheo} because $R_\lam f=R_0(f-\lam R_\lam f)$.
\qed

\begin{remark}\label{rem6}
{\rm If $f\in L^p(\R_+^d,\mu)$, we can extend $f$ to all of $\R^d$ by reflection.
So Corollary \ref{lpcor} also applies when we look at the operators $\sA_i$ and $R_\lam$
operating on functions whose domain is $\R^d_+$.

As with the first order estimates, only minor changes are needed if $R_\lam$ is
the resolvent for $\sum_{i=1}^d a_i \sA_i$, where the $a_i$ are finite positive constants.
}
\end{remark}

\section{Uniqueness}

We need one more estimate, and then  we can complete the proof of
Theorem \ref{maintheorem}.

\begin{lemma}\label{secondderiv}
If $g$ is in $C^2$ with compact support contained in $(0,\infty)^d$, then for each $t>0$  and $\lam>0$ we have
that $P_tR_\lam g$ is $C^2$ in $\R^d_+$ and for each $i$
we have
$(P_tR_\lam g)_i=0$ on $\Delta_i$.
\end{lemma}

\proof We have  a formula for the derivative of the transition density in the
one-dimensional case in
(\ref{ptderiv2}) below in the appendix. If we differentiate once more and use the fact
that the transition densities for the process factor as a product of transition densities
of one-dimensional processes, then
tedious calculations show that $P_tg$ is $C^2$ with normal derivative 0
on the boundary. (This is somewhat easier than in the proof of 
Lemma \ref{l41} since we can use the fact that $g$ has compact support.)
Moreover one can show that the second derivatives of $P_tg$ grow with
$t$ at most polynomially. Since $P_tR_\lam g=\int_0^\infty e^{-\lam s} P_{s+t}g\, ds$
by the semigroup property, the lemma follows.
\qed

The existence part  of Theorem \ref{maintheorem} was done in Section 4. It remains
to prove uniqueness.
\medskip

\noindent{\bf Proof of uniqueness in Theorem \ref{maintheorem}.}
Fix $x_0\in \R^d_+$ and let $\eps>0$ be specified later.
As in the nondegenerate case, to prove uniqueness it suffices to consider only the case
where 
\bee\label{localiz}
 \sum_{i=1}^d |a_i(y)-a_i(x_0)|\leq \eps, \qq y\in \R^d_+;
\eee 
see \cite{Ba97}, Section VI.3.
Let $p_0$ be the positive real given by Theorem \ref{krylovprop}.
Set 
\[
\sL_0 f(x)=\sum_{i=1}^d a_i(x_0)x_i^{\al_i} \frac{\del ^2 f}{\del x_i^2}(x) 
\] and let
\[
\sB=\sL-\sL_0.
\]
Let $R_\lam$ and $P_t$ be the resolvent and semigroup, respectively, for the
operator $\sL_0$. Taking into account Remarks \ref{rem5} and \ref{rem6}, by
Theorem \ref{t46} and Corollary \ref{lpcor} we have
\bee\label{perturb}
\norm{\sB R_\lam f}_{p_0}\leq c_1 \Big(d\eps+\lam^{\frac{1}{2b_i}-1}\sum_{i=1}^d \norm{b_i}_\infty\Big)
\norm{f}_{p_0}.
\eee
Let us now choose $\eps$ small enough and $\lam$ large enough so that by (\ref{perturb})
we have
\bee\label{perturb2}
\norm{\sB R_\lam f}_{p_0}\leq \tfrac12 \norm{f}_{p_0}.
\eee

Let $\P_1$ and $\P_2$ be any two solutions to the submartingale problem
for $\sL$ started at $x_0$, where we continue to assume (\ref{localiz}) holds.
We also assume that under each $\P_i$ the process spends zero time on $\Delta$.
Define
\[ S_\lam ^i h=\E_i \int_0^\infty e^{-\lam t} h(X_t)\, dt, \qq i=1,2. \]
Let $R_K=\inf\{t: \sum_{i=1}^d L_t^{X^i}\geq K\}$.
By Ito's formula, if $f\in C^2$ and $f_i=0$ on $\Delta_i$ for each $i$, then
\[ \E_i f(X_{t\land R_K})-f(x_0)=\E_i\int_0^{t\land R_K} \sL f(X_s) \, ds, \qq i=1,2.
\]
We let $K\to \infty$, so that $R_K\to \infty$. we then multiply both sides by $\lam e^{-\lam t}$
and integrate over $t$ from $0$ to $\infty$ to obtain
\bee\label{sident}
\lam S_\lam^i f-f(x_0)= S_\lam^i \sL f=S_\lam^i \sL_0 f+S_\lam ^i \sB f.
\eee
Now let $g$ be $C^2$ with compact support contained in $(0,\infty)^\infty$ and let $f=P_tR_\lam g$.
By Lemma \ref{secondderiv} we can apply (\ref{sident}) to $f$.
Note
\[ \sL_0 f=\sL_0 R_\lam P_t g=\lam R_\lam P_t g-P_t g.
\]
Therefore (\ref{sident}) becomes
\bee\label{fident}
S_\lam^i P_t g=R_\lam P_t g+ S_\lam ^i \sB R_\lam P_t g.
\eee

Let 
\[ \Theta=\sup_{\norm{h}_{p_0}\leq 1} |S_\lam^1 h-S_\lam^2 h|.
\]
By Theorem \ref{krylovprop} we know $\Theta<\infty$. 
By (\ref{fident}) and (\ref{perturb2}),
\bea
|S_\lam^1 P_t g-S^2_\lam P_t g|&\leq \Theta \norm{\sB R_\lam P_t g}_{p_0}\\
&\leq \tfrac12\Theta\norm{P_tg}_{p_0}\\
&\leq \tfrac12 \Theta \norm{g}_{p_0}.
\end{align*}
The last inequality follows by Jensen's inequality and the fact that $P_t$ is self-adjoint
with respect to $\mu$.
Since the support of $g$ is disjoint from $\Delta$, we can let $t\to 0$
and obtain
\[ |S^1_\lam g-S^2_\lam g|\leq \tfrac12\Theta\norm{g}_{p_0}. \]
We now take the supremum over all such $g$ that in addition satisfy $\norm{g}_{p_0}\leq 1$.
Since neither $S^1_\lam$ nor $S^2_\lam$ charge $\Delta$, we then have
\[ \Theta\leq \tfrac12 \Theta. \]
Since $\Theta<\infty$, we conclude $\Theta=0$.

{}From this point on we follow the proof of the nondegenerate case; see \cite{Ba97}, Section VI.3.
\qed

\section{Appendix}

We prove Lemma \ref{l41}.

We will use the well known facts (see \cite{Lebedev}, pp.~150--152):
\begin{align}
I'_p(x)&= I_{p+1}(x)+\frac{p}{x} I_p(x),\label{Iderivup}\\
I'_p(x)&= I_{p-1}(x)-\frac{p}{x} I_p(x),\label{Iderivdown}\\
I_p(x)&\sim \frac{1}{2^pp!}x^p, \qquad x\to 0, \label{Ismall}\\
I_p(x)&\sim \frac{1}{\sqrt{2\pi}}\frac{e^x}{\sqrt x}, \qquad  x \to \infty. \label{Ilarge}
\end{align}
In what follows we will take $p=\nu$ or $\nu+1$ and $\nu=-1/(2b)$.

If we let $F(x)=I_{\nu+1} (x)-I_\nu(x)$, then from (\ref{Iderivup}) and (\ref{Iderivdown})
we have
\[
F'(x)=-F(x)-\frac{\nu+1}{x} I_{\nu+1}(x)-\frac{\nu}{x} I_\nu(x).
\]
Using (\ref{Ilarge})
\[
|F'(x)+F(x)|\leq c_1 \frac{e^x}{x^{3/2}},
\]
or
\[
|(e^x F(x))'|\leq c_1 \frac{e^{2x}}{x^{3/2}}.
\]
Therefore
\[
|e^x F(x)|\leq |e F(1)| +c_1\int_1^x \frac{e^{2y}}{y^{3/2}} dy.
\]
By l'H\^opital's rule, the integral is bounded by $c_2 e^{2x}/x^{3/2}$, and so
we deduce
\begin{equation} \label{Idifference}
|I_{\nu+1}(x)-I_\nu(x)|\leq c_3 \frac{e^x}{x^{3/2}}
\end{equation}
for $x\geq 1$.

\medskip
\noindent{\bf Proof of Lemma \ref{l41}.}  We start with (\ref{ap23}).
By scaling it suffices to do the case $t=1$. Differentiating (\ref{ptdef}) we
have
\bee\label{ptderiv2}
\frac{\del p_X(1,x,y)}{\del x}=cx^{b-\frac12}y^{2b-\frac32} e^{-K(x^{2b}+y^{2b})/2}
[-x^b I_\nu(Kx^by^b)+y^bI_{\nu+1}(Kx^by^b)],
\eee
where $K$ is some fixed positive constant.
We write
\[
\int_0^\infty \Bigl|\frac{\del p_X(1,x,y)}{\del x}\Bigr|\,dy=\int_0^{1/x}+\int_{1/x}^\infty:=S_1+S_2.
\]
Using the bounds on $I_\nu$,
\[
S_1\leq cx^{2b-1}e^{-Kx^{2b}/2}\int_0^{1/x} [y^{2b-2}+y^{4b-2}]e^{-Ky^{2b}/2} dy.
\]
Since $2b-2>-1$, the integral term is finite. Since $2b-1>0$, the factor in front
of the integral is bounded independently of $x$, so $S_1$ is bounded independently
of $x$.

Since
\begin{align*}
|-x^bI_\nu(Kx^by^b)&+ y^bI_{\nu+1}(Kx^by^b)|\\
&= |(y^b-x^b) I_\nu(Kx^by^b)+y^b(I_{\nu+1}(Kx^by^b)-I_\nu(Kx^by^b))|\\
&\leq c|y^b-x^b| e^{Kx^by^b} x^{-\frac{b}{2}} y^{-\frac{b}{2}}+
y^b(e^{Kx^by^b}x^{-b}y^{-b})
\end{align*}
for $y\geq 1/x$, to bound $S_2$ we need to bound
\begin{align*}
\int_{1/x}^\infty& |y^b-x^b| x^{\frac{b-1}{2}} y^{\frac{3b-3}{2}} e^{-K(y^b-x^b)^2/2} dy\\
&~~~+\int_{1/x}^\infty x^{-\frac{1}{2}} y^{2b-\frac{3}{2}} e^{-K(y^b-x^b)^2/2} dy\\
&= S_3+S_4.
\end{align*}
For $S_3$ we make the substitution $z=y^b-x^b$ and then
\[
S_3=c\int_{x^{-b}-x^b}^\infty |z| x^{\frac{b-1}{2}}(x^b+z)^{\frac{b-1}{2b}} c^{-kz^2/2} dz.
\]
Since $(b-1)/(2b)<0$ and $x^b+z\geq x^b$, this is less than
\[
c\int_{-\infty}^\infty |z| x^{\frac{b-1}{2}} (x^{-b})^{\frac{b-1}{2b}} e^{-Kz^2/2} dz
\]
which is bounded independently of $x$.  For $S_4$ we use the same substitution. Since $2b-1>0$,
we have 
\[
(x^b+z)^{\frac{2b-1}{2b}}\leq c(x^{\frac{2b-1}{2}}+z^{\frac{2b-1}{2b}}).
\]
Hence 
\begin{align*}
S_4&\leq  \int_{x^{-b}-x^b} ^\infty x^{-\frac12} (x^b+z)^{\frac{2b-1}{2b}} e^{-Kz^2/2}dz\\
&\leq   c\int_{x^{-b}-x^b}^\infty 
(x^{b-1}+x^{-\frac12}z^{\frac{2b-1}{2b}})e^{-Kz^2/2} dz.
\end{align*}
For each $p\geq 1$ there exists $c(p)$ such that
\bee\label{tail}
\int_a^\infty (1+z)^{\frac{2b-1}{2b}} e^{-Kz^2/2} dz\leq c(p) a^{-p}, \qquad a>1.
\eee
{}From this we see that $S_4$ is bounded independently of $x$ for $x\leq1$. On the other hand,
for $x\geq 1$,
\[
S_4\leq c\int_{-\infty}^\infty (1+|z|^{\frac{2b-1}{2b}}) e^{-Kz^2/2} dz\leq c.
\]

We now turn to the proof of (\ref{ap25}).
Again by scaling we may assume $t=1$. 
Looking at $\int_0^{1/y}$ and  using the bounds on $I_\nu$, 
\begin{align*}
 y^\al &\int_0^{1/y} \Bigl|\frac{\del p_X(t,x,y)}{\del x}\Bigr| x^{-\al} dx\\
&\leq  c[(y^{\al+2b-2}+y^{\al+4b-2}) e^{-Ky^{2b}/2}] \int_0^{1/y}
x^{2b-1-\al} e^{-Kx^{2b}/2} dx.
\end{align*}
Since $\al+2b-2=0$, $\al+4b-2>0$, and $2b-1-\al=1-2\al>-1$, the
integral is finite and the expression in brackets is  bounded in $y$.

To look at $\int_{1/y}^\infty$, we rewrite the integral as in $S_2$
and see that we have to bound
\begin{align*}
cy^{\frac32 b-\frac32+\al}&   \int_{1/y}^\infty
x^{\frac{b}2 -\frac12 -\al} |y^b-x^b| e^{-K(y^{b}-x^{b})^2/2} dx\\
  & \qquad +cy^{2b-\frac32+\al}\int_{1/y}^\infty x^{-\frac12-\al}
e^{-K(y^{b}-x^{b})^2/2} dx\\
&=  S_5+S_6.
\end{align*}
Letting $z=x^b-y^b$ as in $S_3$,
\begin{align*}
S_5&\leq   cy^{\frac32 b-\frac32+\al}\int_{y^{-b}-y^b}^\infty
(y^b+z)^{-\frac{b+1+2\al}{2b}} z e^{-z^2/2} dx\\
&\leq   c y^{b-1} \int_{y^{-b}-y^b} z e^{-z^2/2} dz.
\end{align*}
When $y\leq 1$ this is bounded using (\ref{tail}). When $y\geq 1$
this is bounded because $b-1<0$.

For $S_6$ we have
\begin{align*}
S_6&\leq   cy^{2b-\frac32+\al} \int_{y^{-b}-y^b}^\infty (y^b+z)^{-\frac{b-3-2\al}{2b}} e^{-Kz^2/2} dz\\
&\leq  c y^{\frac{b}2 -3} \int_{y^{-b}-v^b}^\infty e^{-Kz^2/2} dz.
\end{align*}
This is bounded in $y$ for $y\leq1$ by (\ref{tail}). This is bounded in $y$
for $y>1$ because $\frac{b}2-3$ is negative.
\qed

\def\polhk#1{\setbox0=\hbox{#1}{\ooalign{\hidewidth
  \lower1.5ex\hbox{`}\hidewidth\crcr\unhbox0}}}
  \def\polhk#1{\setbox0=\hbox{#1}{\ooalign{\hidewidth
  \lower1.5ex\hbox{`}\hidewidth\crcr\unhbox0}}}
  \def\polhk#1{\setbox0=\hbox{#1}{\ooalign{\hidewidth
  \lower1.5ex\hbox{`}\hidewidth\crcr\unhbox0}}} \def\cprime{$'$}
  \def\cprime{$'$}


\medskip

\begin{minipage}[t]{0.39\textwidth}
{\bf Richard F. Bass}\\
Department of Mathematics\\
University of
Connecticut \\
Storrs, CT 06269-3009, USA\\
{\it bass@math.uconn.edu}
\end{minipage}
\hfill
\begin{minipage}[t]{0.55\textwidth}
{\bf Alexander Lavrentiev}\\
Department of Mathematics\\
University of Wisconsin, Fox Valley\\
Menasha, Wisconsin 54952, USA\\
{\it alavrent@uwc.edu} 
\end{minipage}

\end{document}